# The Technological Turn in Mathematics

Silvia De Toffoli and Fenner Tanswell

ABSTRACT: Quickly evolving technologies, such as Interactive Theorem Provers (ITPs), Automated Theorem Provers (ATPs), and Large Language Models (LLMs), all falling under the general heading 'AI for mathematics,' are transforming mathematical practice in profound ways. This chapter explores the implications of these innovations, focusing on their impact on how mathematical knowledge is created and shared. It also discusses how they are reshaping the social dimension of mathematics, altering collaboration dynamics, trust relationships, and the collective production of knowledge. For instance, tools like ITPs facilitate large-scale collaborations and make new types of teamwork possible, where trust is not a necessary ingredient. ITPs also help us mitigate our human fallibility, yet they raise questions about the nature of formalization and the relationship between traditional and formal mathematics. Technologies such as LLMs are reshaping the division of epistemic labour between humans and machines and urge philosophers of mathematics to ask questions about the value of their work.

## 1. Introduction

New technologies are increasingly affecting mathematical practice, radically transforming how new results are discovered and justified. This transitional period can be characterized as a genuine "technological turn," significantly altering the epistemology of mathematics as well as its social structure.

One might think that this is an overstatement. After all, the use of different technological tools has been integral to mathematical research for a very long time. The abacus was already used throughout the ancient world. Mechanical tools designed as calculation devices such as Pascal's calculator and Leibniz's stepped reckoner, were already developed in the seventeenth century. Digital calculators have been widely available to mathematicians since at least the 1970s. Computer Algebra Systems (CAS) such as *Mathematica* and *Maple* became widely used in the 1980s, along with numerical computing environments such as *MATLAB*. So, the use of digital technological tools is by no means a novelty in mathematics.

Nothing new even within the heart of research mathematics, it might seem. Already in the late 1970s, Appel and Haken (1977) devised what is often considered to be the first computer-assisted proof.[1] They proved the Four Color Theorem with the aid of a computer. And in the 1990s, Doron Zeilberger (1994) claimed that "[t]he computer has already started doing to mathematics what the telescope and microscope did to astronomy and biology" (11).[2] However, it is only in recent years that the use of Artificial intelligence in mathematics (and in other domains) has exploded. In another chapter of this volume, Tanya Klowden and Terence Tao (forthcoming) share the feeling that this is a turning point: "Humanity is standing at the threshold of a digital Industrial Revolution, unfolding at unprecedented speeds" (1).

'Artificial Intelligence in mathematics' is an umbrella term encompassing a wide range of different tools. Adapting Jeremy Avigad's (2025) taxonomy,[3] we can identify three main clusters of technologies that fall under it:

a) Interactive Theorem Provers (ITPs) used for formalizing mathematics,
b) Symbolic AI used for automating mathematical reasoning, and
c) Machine learning and neural AI.

In this chapter, we mostly focus on philosophical questions raised by the introduction of a) and c) within mathematical practice since they are the technologies that have advanced more dramatically in recent years. We consider this the titular 'technological turn' in mathematics, and will consider how it affects a range of philosophical issues on proof, justification, rigor, fallibility, collaboration, trust, and automation.

Our paper proceeds as follows. We begin in Section 2 with a brief overview of the relevant technologies. Then we will deal with the following emerging philosophical issues: Section 3 deals with how ITPs are changing what mathematical justification amounts to. It discusses the nature of proof, mathematical rigor, and the prospect of a fallibilist epistemology of mathematics in the age of ITPs. Section 4 continues the analysis of ITPs, this time from the perspective of the new types of collaborations they enable. Section 5 considers how ITPs might alter the structure of the trust relationships in mathematics. Section 6 zooms in on the use of neural AI in mathematics, focusing on both the hopes and the concerns it generates among professionals. Section 7 concludes with some final considerations.

---

[1] Although, to be clear, there were numerous earlier proofs involving computers, as described in (Detlefsen and Luker 1980).

[2] Zeilberger argued that computers would lead mathematicians to abandon the idea of perfect rigor – this did not exactly go as he thought. For a philosophical discussion, see (De Toffoli 2024).

[3] See also Avigad's (forthcoming) chapter in the present volume.

## 2. A Brief Overview of the Technology

Older technologies are now used alongside newer ones arising from rapid advancements in both symbolic and neural AI. These two types of AI model intelligence in radically different ways. In a sentence, the paradigm behind the first is logic-based, while that behind the latter is brain-inspired (Yoshua Bengio, LeCun, and Hinton 2021). In more detail, symbolic AI, sometimes called 'Good Old-Fashioned AI' (GOFAI), grew out of the mathematical logic tradition and its influence on early computer science. This paradigm for AI is characterised by mathematical symbolic manipulations, representative abstractions, compositionality, logical inferences, and explicit knowledge representations in databases, maps, semantic networks, and libraries. Neural AI is instead based on *neural networks*, loosely inspired by the structure of the human brain. This paradigm is characterised by statistical inference, training on large amounts of data, reinforcement learning, uncertainty and unpredictability, and implicit knowledge encoded in weighted networks. The symbolic paradigm underpins most modern programming languages and software while the neural paradigm underlies the transformer architecture used by large language models (LLMs).

The symbolic tradition has numerous milestone successes, such as Heule, Kullmann, and Marek's (2016) 200-terabyte solution to the Boolean Pythagorean triples problem.[4] However, the modern surge in interest in AI in mathematics is generally driven by advances in LLMs and deep learning. In particular, the development of neural AI has witnessed a breakthrough since the popularization of large language models (LLMs) such as ChatGPT at the end of 2022. AI systems of this type were first designed to work with natural language and did not do well with mathematics, but things are rapidly changing as we will see below.

It should be noted that there is the possibility (and hope) of combining the two approaches with a *neurosymbolic* approach, combining the advantages of neural AI with the rigorous reasoning that underlies symbolic AI. For example, DeepMind's *AlphaGeometry* model is described by Trinh et al. (2024) as "a neuro-symbolic system that uses a neural language model, trained from scratch on our large-scale synthetic data, to guide a symbolic deduction engine through infinite branching points in challenging problems" (Trinh et al. 2024, 476).

[4] This is the question of whether the positive integers can be divided into two sets so that no Pythagorean triples are from one set, where a Pythagorean triple consists of three positive integers $(a, b, c)$ such that $a^2 + b^2 = c^2$. Heule, Kullmann, and Marek's (2016) solution shows that this is possible up to 7824 but no higher.

## 3. Proof, Rigor, and Fallibility

Mathematical rigor is often taken to be the essence of mathematical method: proofs must be rigorous. According to the Standard View, genuine proofs are arguments that can be formalized in an appropriate formal system (Burgess and De Toffoli 2022; Hamami 2022; Tanswell 2024). In one form or another (and indeed it admits many forms), this view is widely accepted among mathematicians. If a purported proof were shown to be impossible to formalize (in an appropriate formal system), that would generally be taken to be evidence that it was not a proof after all.

Formalizability is thus often taken to just mean correctness (Avigad 2021). However, formalization is not generally used to check for correctness (De Toffoli 2021b) as humans are not very good at verifying whether an argument is formally correct, and much better at verifying traditional mathematical arguments, which are informal arguments written in a combination of natural language and different notational systems. And in fact, mathematicians generally check the correctness of their (informal) proofs by thinking through them.

Moreover, mathematicians often rely on one another to verify the correctness of their results. Although a proof guarantees that the conclusion follows from the premises, mathematicians know they are not infallible in recognising proofs. They are humans and like all other humans make mistakes. This is why there are cases in which a mathematician announces a proof without having one. Sometimes an incorrect argument containing an invalidating error is mistaken for a proof. The error could be either obvious, and thus immediately spotted when the putative proof is subjected to external scrutiny, or very subtle, and elude other experts as well, or something in between.

It is plausible to think that linking mathematical justification to proof alone is too strong a requirement – note that "proof" is a success term and thus there is no such thing as an "erroneous proof," even if sometimes this expression is used in practice.[5] According to De Toffoli (2021a, forthcoming-a, forthcoming-b), in order to be justified in believing a new mathematical result, a mathematician needs to have a "simil-proof," that is, an argument that conforms to the inferential standards of a legitimate mathematical community (and thus generally is taken to be a proof by the relevant experts). Simil-proofs are often genuine proofs, but sometimes they are erroneous putative proofs containing subtle mistakes. In practice, the term "proof" is used loosely to mean "simil-proof." For instance, when a mathematician initially publishes a proof, at first glance we can only gain evidence that her argument is a simil-proof (see Section 5). Often, this simil-proof will turn out to be a genuine proof, but in some cases, it will turn out to contain a subtle flaw.

[5] See (Paseau 2015) for a defence of the claim that in some cases it is possible to achieve mathematical knowledge without proof.

The first simil-proof of the Four Color Conjecture (which many called a “proof”) was not a proof but an erroneous simil-proof. It was published by Alfred B. Kempe in 1879 and broadly accepted by the mathematical community for eleven years, when a subtle gap was found. As Kempe’s case suggests, community checking is crucial to filtering out erroneous simil-proofs, ensuring the reliability of the mathematical literature. However, we now face a new challenge. Mathematics is becoming increasingly complex, and simil-proofs, especially involving automated computations, are getting harder to check:

> With results like the four color theorem and Hales’ theorem, we are gradually getting past the vain hope that every interesting mathematical theorem will have a humanly surveyable [simil-]proof. But it seems equally futile to hope that every computational [simil-]proof will make use of code that can easily be understood, and so the usual difficulties associated with understanding complicated proofs will be paired with similar difficulties in understanding complicated programs. (Avigad 2018, 687)

Even without automated computation, simil-proofs in some mathematical fields are so long and complex that often no one checks them. This is a problem because there might be errors in accepted results, and if these errors are not corrected, they tend to generate new errors, affecting entire areas of mathematics. How can we make sure that the simil-proofs we accept are actually proofs and do not contain subtle errors?

Mathematicians are increasingly relying on *Interactive Theorem Provers* (also called ‘proof assistants’) to check the correctness of their proofs. These are software programs that allow mathematicians to convert their informal arguments, written in natural language and mathematical notation, into formal arguments in a specific formal system (Buzzard 2024; Hales 2008). These formal arguments can then be mechanically checked for correctness. Similar tools are also used to formally verify software, used both within and outside mathematics. In the words of one of the pioneers of formalization, George Gonthier (2008),

> The idea is to write code that describes not only *what* the machine should do, but also *why* it should be doing it—a formal proof of correctness. The validity of the proof is an objective mathematical fact that can be checked by a different program, whose own validity can be ascertained empirically because it runs on many inputs. (1382)

ITPs are “interactive” because it is mathematicians who construct the proofs, choosing strategies, supplying key ideas, etc., while the system ensures that every step follows from precise logical rules. Examples of such systems are *Lean*, *Coq* (and its recent development, *Rocq Prover*), and *Isabelle*. One key difference among these is that they adopt different foundational systems: for instance, Coq and Lean are based on dependent type theory, whereas Isabelle is primarily based on higher-order logic.[6]

One of the first high-profile formalization projects using an ITP was the formalization of Appel and Haken's proof of the Four Color Theorem mentioned above. After all, the early example of a computer-assisted proof had raised suspicion among mathematicians and philosophers alike about its status and validity (Tymoczko 1979). One of the problems was that, not being surveyable, it was impossible for human agents to check all the details. That is why this proof was a perfect candidate for formal verification. Such verification was achieved by Gonthier in 2005 using Coq (Gonthier 2008).

One early proponent of these tools was Fields Medalist Vladimir Voevodsky. Working in a field where "'too-long [simil-]proofs' are common and in which the relatively small number of competent potential referees typically spend much of their time writing 'too-long [simil-]proofs' of their own" (Harris 2015, 58), he started to worry that the self-monitoring activity of the community was not in place anymore. Another hint that something was going awry with community checking was that a counterexample to one of his simil-proofs was put forward, but nobody could find a flaw in his simil-proof until a decade later. That is why he turned to ITPs as an alternative means of ensuring the reliability of the mathematics literature (Voevodsky 2014).

Another early breakthrough is Thomas Hales' verification of his proof of the Kepler conjecture. Hales announced that he had proven it in 1998. His proof was published some years later in the *Annals of Mathematics* (Hales 2005). However, the status of his argument remained unclear because the referees found it difficult to check all the details and only accepted it with some reservations because, like Appel and Haken's proof, Hales's proof was also computer-assisted. Motivated by the desire to establish the validity of his proof, he and his collaborators began the enterprise of formally verifying it, which they completed in 2017 (Hales et al. 2017).

Driven by similar concerns, Peter Scholze, another Fields Medalist, appealed to the Lean community in 2022 to check the correctness of one of his simil-proofs.[7] As it turned out, no one else checked all the details (Buzzard 2021, 12). Scholze's call was enthusiastically answered, and his doubts were surprisingly quickly put to rest by a team led by Johan Commelin. A high-profile project currently underway led by Kevin Buzzard aims to provide a full formalization of Andrew Wiles and Richard Taylor's proof of Fermat's Last Theorem.[8]

---

[6] For a survey, see (Avigad 2018).

[7] The proof was a foundational result within his condensed mathematics program (Scholze 2019). The Lean community responded enthusiastically and started the *Liquid Tensor Experiment*, which was extremely successful. See: https://xenaproject.wordpress.com/2020/12/05/liquid-tensor-experiment/ and https://leanprover-community.github.io/blog/posts/lte-final/.

[8] See: https://leanprover-community.github.io/blog/posts/FLT-announcement/.

ITPs can thus ensure the correctness of our simil-proofs when community checking is missing. Moreover, even when community checking is in place, ITPs provide more thorough scrutiny. As individual mathematicians at times err, groups of mathematicians can also make mistakes (e.g. Kempe's putative proof of the Four Color Conjecture was accepted by the mathematical community for 11 years before the problem was found). Full formalization achieved through ITPs can help eliminate both individual and community errors.

According to some, in the near future, at least in some areas of mathematics, regular simil-proofs (especially if computer-assisted in some way) will have to be paired with a formal counterpart (Avigad 2024; Yang et al. 2024). This might lead to a change in the norms governing mathematical practice. Some mathematicians are already starting to formalize their new simil-proofs. For example, after the publication of a preprint containing a simil-proof of the Polynomial Freiman-Ruzsa (PFR) conjecture over $\mathbb{F}_2$ (Gowers et al. 2024), Terence Tao, one of the authors of the simil-proof, started a collaborative project to formalize it.[9]

As we saw, ITPs can help us mitigate our human fallibility even when social mechanisms of self-monitoring are not feasible. Does this mean that we can now rest assured and be completely certain of the correctness of our results? Although this might be an enticing idea, there are reasons to be more cautious as new forms of *technological fallibility* become relevant.

ITPs are software, and like all software, might contain bugs. Critical bugs are routinely found and fixed in many proof assistants.[10] As a matter of fact, someone has even "proved" an inconsistent statement to show that the system contained a bug:

> A couple of days ago, Daniel Litt linked to Patrick Brosnan's computer-verified "proof" of the inconsistency of Peano arithmetic. The proof is correct; I just put it in quotes because it relies on a quirk of the proof verification system used (Metamath), which requires you to explicitly prohibit certain variable substitutions. (Zach 2021)

A subtler issue has to do with the relationship between formal and informal mathematics. How can we be entirely sure that what we mean to formalize actually corresponds to what we do in fact formalize? This is a classic problem in the philosophy of mathematical practice, considered at length by Lakatos (1976) as the problem of translation.[11] It is also a practical concern for mathematicians involved in formalization efforts. In another chapter in this volume, Klowden and Tao (forthcoming) write:

> [F]ormal verification only certifies that a formalized argument establishes a formal mathematical statement, but does not rule out errors in translation between the formal statement and the original intended statement. (9)

[9] https://terrytao.wordpress.com/2023/11/18/formalizing-the-proof-of-pfr-in-lean4-using-blueprint-a-short-tour/.

[10] See, for example, https://github.com/rocq-prover/rocq/blob/master/dev/doc/critical-bugs.md, which "recollects knowledge about critical bugs found in Coq since version 8.0."

For instance, in certain mathematical domains, such as differential geometry, it is non-trivial to define the main objects of study in their full generality (e.g., differentiable, smooth, analytic, or complex manifolds, or manifolds endowed with additional structure such as Lie groups). It is therefore important to ensure that no substantive assumptions are silently built into the very definitions being formalized. The risk is not merely that a formalization fails, but that it inadvertently establishes a different result than intended. This problem can be exacerbated when interactive theorem provers are paired with neural AI systems, such as Aristotle. When formalized statements are passed to AI systems trained to search for proofs, these systems may exploit hidden assumptions in the formal setup, yielding results that are formally correct but conceptually misaligned, in that they use a definition other than the intended one.[12]

To sum up, although successful formalization via an ITP does justify greater confidence in the correctness of a simil-proof, it may still fall short of absolute certainty.[13]

## 4. Collaboration, Crowdsourcing, and Hybrid Systems

In section 3, we discussed some of the epistemic changes ITPs are bringing to mathematical practice. Now, we turn to ITPs' social consequences, focusing on the fact that they enable carrying out mathematical collaboration at much larger scales.

While there have been examples of large collaborations without substantial dependence on ITPs, most notably the enormous proof of the Classification of Finite Simple Groups (Habgood-Coote and Tanswell 2023; Steingart 2012), these have the problem that it is inevitable that errors will arise in such a large (simil-)proof. The general opinion of the group theory community is that all of the remaining errors in the proof of the Classification Theorem are merely "local errors" (as discussed by Steingart (2012)), which are easily fixable and do not threaten the correctness of the overall proof. Interestingly, the confidence that none of these supposedly local errors hides a major problem is not based on the argument itself, but on the thought that if there were a major (non-local, invalidating) error, somebody would have noticed it.[14]

[11] See also (Burgess and De Toffoli 2022; Tanswell 2015; 2024, Sec. 3.7).
[12] Thanks to Tito Sacchi for helping us make this point more precise.
[13] See also the discussion in (De Toffoli 2024).

Still, there does not appear to be any guarantee that there are no errors that are both major and *subtle*: those that would not have been spotted even by repeated checking by different mathematicians. That is, there is no guarantee that this is not an erroneous simil-proof. To be sure, the problem arises even with arguments produced by a single mathematician, but in very large proofs that require collaboration, it is exacerbated. Furthermore, a single mathematician can be reasonably expected to maintain an overview of how the entire putative proof coheres, something which may no longer be the case for a large collaboration.

ITPs provide the infrastructural support to enable modern large-scale collaborations, while potentially avoiding worries about large-scale mathematics being swamped by errors. By using ITPs, every contribution must be written up to the demanding standards of the relevant formal language and thus can be formally verified. Here are four benefits ITPs bring to large-scale collaborations.

First, trust is no longer essential, as emphasized in (Tao 2025b). By formally verifying results, mathematicians no longer have to trust that the work of their collaborators is sufficiently error-free to use in their own work. This has the potential of democratizing mathematics since it makes collaborations possible where no antecedent personal trust is in place, so also no need to rely on a mathematician's reputation or other imperfect social markers. Second, the overall structure is available and accessible. The library software can track mathematical dependencies and connections, allowing mathematical work to build cumulatively without any single mathematician needing to have it all in mind at once. Moreover, the library is searchable, making it much easier to find lemmas and intermediate results. Third, the systems are widely accessible and stable, so collaboration can happen asynchronously with researchers around the world at any time. Fourth, the systems are such that mathematics can be modularised, and therefore distributed to different mathematicians to work on, then recombined without the danger of some mismatch in their work causing an error.

ITPs thus have the potential to reshape the social organisation of mathematics by allowing for smoother *crowdsourcing* of mathematical work. According to a systematic overview of definitions, crowdsourcing refers to "a type of participative online activity in which an individual, an institution, a non-profit organization, or company proposes to a group of individuals of varying knowledge, heterogeneity, and number, via a flexible open call, the voluntary undertaking of a task" (Estellés-Arolas and González-Ladrón-de-Guevara 2012, 197). Pease et al. (2019) consider various examples of mathematical crowdsourcing, such as asking questions on mathematics-themed question and answer sites like *MathOverflow*, *Math.Stackexchange*, or reddit's */r/math*. Mathematical crowdsourcing does not necessarily need to involve the use of ITPs, and indeed, early examples did not include them. However, given the benefit they provide, ITPs can bring mathematical crowdsourcing to the next level.

[14] That is, this would be what Goldberg (2011) calls a *coverage-supported* belief. See the discussion in (Habgood-Coote and Tanswell 2023).

Let us now consider two examples of mathematical crowdsourcing, the *Polymath Projects*, which did not involve ITPs, and the more recent *Equational Theories Project,* which includes the use of an ITP. The *Polymath* Projects sought to make progress on a series of research-level mathematics questions through crowdsourcing.[15] The project started on the blog of the Cambridge mathematician Timothy Gowers in 2009, where he proposed an open, online, collaborative effort to give a new, elementary proof of the Density Hales-Jewett theorem (Gowers 2010; Gowers and Nielsen 2009). The Polymath Projects were run as online comment threads and structured around this medium (Barany 2010). In Gowers's original set of guidelines for the project, he encouraged all work to be put online and done publicly, instead of individuals reporting on work they have done in private. The experiment in crowdsourced mathematics was a success, both resulting in a successful proof (Polymath 2012) and extensive online collaboration:

> Over the next seven weeks, 27 people contributed around 800 comments- around 170,000 words in all- with the contributors ranging from high-school teacher Jason Dyer to Gowers's fellow Fields Medallist Terry Tao. (Martin 2015, 32)

By now there have been sixteen Polymath Projects, with the last one, on making progress on the Hadwiger–Nelson problem, finished in 2021. Not all of the projects were successful, with some winding down due to a lack of progress or engagement, but numerous new mathematical results were achieved.

In contrast to the linear blog style of the early Polymath Projects, a 2024 project allowed for a different, less linear, division of labour. In September 2024 Terence Tao, who had been frequently involved in the Polymath Projects, proposed a crowdsourced project he called the *Equational Theories Project* (ETP).[16] The project was about *magmas*: a set $M$ with a binary operation $\diamond\colon M \times M \to M$. Equational laws are then laws stated using the operation, e.g., *Commutativity* is the law $x \diamond y = y \diamond x$, and *Tarski's axiom* is the law $x = y \diamond \big(z \diamond \big(x \diamond (y \diamond z)\big)\big)$. Different equational laws can be combined together to form equational theories, but only some of these will be consistent, i.e., some combinations of equational laws will imply that others are false. The aim of the project was to determine the implication graph for all of the equational theories that use the 4694 equational laws that used the binary operation up to 4 times, which leaves an implication graph with 22,028,942 edges.

[15] See (Pease et al. 2020).

[16] Our description of the project follows the paper produced by the participants at the end of the project (Bolan et al. draft). We focus on this case here, but there are several similar projects, such as the *Busy Beaver Challenge,* and Kevin Buzzard's *Fermat project* crowdsourcing of parts of the Lean formalisation of Fermat's Last Theorem (see Section 4).

The idea was to provide a mathematical challenge that could be properly crowdsourced and distributed among participants who could work on different implications in parallel. Furthermore, it was clear that large parts of the work should be automated. Rather than the blog style, which requires moderation and constant overview, as well as unscalable human integration of results, the outputs were all formalised in Lean and formally verified. Outputs were not added to a blog but instead to a *GitHub* repository.

The full implication graph was finished in April 2025.[17] The final project report lists thirty contributors as authors, and states that fifty were involved overall. They describe the extensive use of different kinds of software: Automated Theorem Provers such as *Vampire, Mace, Duper, Z3, MagmaEgg,* and *Prove9,* were used to automate large parts of the mathematics, Lean was used for verification, Patrick Massot's *Lean blueprint* system was used for project management, a dedicated section of the *Lean Zulip* messaging forum, *GitHub's* copilot feature for code autocomplete suggestions, and various pieces of custom interface software including the *ETP dashboard*, the *Equation Explorer*, *Graphiti*, and the *Finite Magma Explorer*.

The difference between the Equational Theories Project and earlier Polymath Projects is that the social epistemic picture that emerges is much more complex. In linear blogpost mathematics, the social epistemology is one of large-scale collaboration, but where the final results are generally surveyable by all participants. That means that the knowledge may be *cumulative* in the sense discussed in (Habgood-Coote and Tanswell 2023): "between them the group theory community has objectual knowledge of the complete proof, by its members having knowledge of parts of that proof" (17), but in practice this may be because some individual in the group knows the entire proof.

In the Equational Theories Project, however, the situation is better conceptualised as a *hybrid epistemic network*, combining both human mathematicians, automated computer tools, and technological infrastructure, in order to complete a unified mathematical project. This begins to realise Ursula Martin and Alison Pease's vision of mathematics as a social machine:

> [A] combination of people, computers, and mathematical archives to create and apply mathematics, with the potential to change the way people do mathematics,

[17] Interestingly, the project also looked at what happened when restricting to finite magmas. Here there were precisely two implications that resisted proof or refutation with any techniques developed in the rest of the project.

and to transform the reach, pace, and impact of mathematics research. (Martin and Pease 2013, 1)

This new way of collaborating with people and machines raises a number of philosophical problems, as it combines aspects of the traditional cases of large collaboration with questions about computer mathematics, as have been discussed extensively in the philosophy of mathematics literature, originally in relation to the computer-aided proof of the Four Color Theorem.

The situation of a hybrid epistemic network, or a social machine of mathematics, is more complicated. The work on the proof has been distributed between human participants, multiple different pieces of code in multiple languages, Lean, and Lean's library *Mathlib*. These parts can be considered the *nodes* in this large epistemic network. However, this model suggests something more discrete than the case we have examined. The real practice involved significant interaction between the different nodes: between human participants, between humans and computers, and between different computer systems. As such, rather than all of the parts of the proof being known by distinct nodes, plenty of parts are only known by a subsection of the network, i.e., by some human-computer hybrid agent.[18]

One social epistemological option is to ascribe mathematical knowledge to the group or collective of mathematicians.[19] In this case, though, the collective is the hybrid epistemic network, or social machine of mathematics, and explicitly includes non-human components. Whether this is an acceptable solution for an epistemology of mathematics is clearly up for philosophical debate.

## 5. Trust and Verification

Another central philosophical issue in large-scale technology-assisted mathematics concerns the role of trust in professional mathematical contexts. In our previous work (De Toffoli and Tanswell 2025), we defended an account of trust in mathematics starting with a critique of the *trust-free ideal of mathematics*. This ideal is that there should be no trust in mathematics, such that every mathematician should check every theorem or result they rely on in their work themselves. While in theory this is no bad thing, we argue that in practice this is limiting mathematics for two main reasons.

[18]One question is what exactly the knowledge is here. It may be that this is purely propositional knowledge of mathematical facts. However, it may also be more: Tanswell (2024) defends the recipe model of proofs, where proofs are the recipes for a series of proving activities. See also (Tanswell and Inglis 2024; Weber and Tanswell 2022). The epistemology of this involves a mathematician needing a combination of propositional knowledge (knowledge *that*) and practical knowledge (knowledge *how*). To gain *a priori* knowledge of a theorem, mathematicians need to know how to carry out the required actions in a proof and then enact them. In this case, the knowledge of how to carry out the proof, and the actual doing, are widely distributed amongst the hybrid network, and no individual could know how to do all of the parts of the projects.

[19] See the discussion in (Habgood-Coote and Tanswell 2023). Similar solutions have been proposed in general epistemology as well. In particular there is a debate on the nature of group knowledge and whether this is reducible or not to individual knowledge (see, for example, (Brown 2024; Lackey 2020)).

First, the division of epistemic labour is essential at larger scales of mathematics, and for this the trust-free ideal is a major hindrance:

> [C]ontemporary mathematics is huge in terms of the quantity of results, definitions, and techniques that have been produced. Mathematics is divided into a myriad of fields and subfields, and communication between them is not always easy. It is far too large for any individual mathematician to have any hope of learning even a substantial fraction of it. Furthermore, mathematics is deeply interconnected, meaning that results from different subfields are often applied to other subfields. Especially in some areas of mathematics, this often makes it impossible (or highly impractical) for mathematicians to check all the results that they rely on. (De Toffoli and Tanswell 2025, 5)

Second, the current incentive structure of mathematics rewards new mathematics and thus makes it extremely difficult for mathematicians (especially junior scholars who have to find a stable location in the profession) to adhere to the no-trust ideal. Moreover, the risk of seemingly well-established mathematics being wrong is often low compared to the huge time and effort required to check it for oneself, time which could be better spent on new research. An evaluation of the optimal course of action clearly favours trusting others over verifying everything oneself (at least in general).

Accepting these considerations, the next question is what exactly trust amounts to in mathematics. We developed an account of trust in mathematics drawing on Katherine Hawley's *commitment account of trust* (Hawley 2014). According to her account, to trust someone is to rely on them to fulfil a commitment. This allows for a rich analysis of the commitments that mathematicians take on at various stages of their practices. In our analysis (De Toffoli and Tanswell 2025), we argue that when a mathematician claims to have a novel proof, this cannot commit them to *actually* having a genuine proof because humans are fallible and sometimes mistakes are subtle enough to be excusable. Rather, the mathematician makes a commitment that their proposed proof meets the professional standards of the discipline without any obvious errors, that is, that their purported proof is a simil-proof (De Toffoli 2021a, 2026). Furthermore, it commits them to having the proof validated by their peers, highlighting the importance of peer review and post-review checking in mathematics, something Helen Longino (1990) famously stressed for science. Conversely, the mathematical community also has a commitment to engage with the work.

In the case of large-scale, technology-augmented mathematical collaboration, the role of trust appears to be quite different. There is no longer the need to trust your collaborators not to have made a mistake because the agreed standard of correctness is that the proofs are formalised and checked by an ITP. No longer is there a danger that there is a subtle error in the work of a collaborator that you need to rely on, since these errors will not make it past the computer verification (at least up to the technological fallibility mentioned in Section 3). There is, therefore, also no apparent need to rely on checking by peers, nor does there need to be a commitment from the community to do so, so also no need to trust that others will engage with your work.

We might wonder whether this setup requires that mathematicians instead trust in the computer. However, computers do not make commitments, so while they are relied upon in the process of doing mathematics, this kind of mere reliance is not the same as trust. Indeed, this account of trust is restricted explicitly to interpersonal trust, so it excludes computers from the outset. On other accounts of trust, though, this is not the case. For example, Coliva (2025) develops a general account of "hinge-trust," where trust is not only interpersonal,[20] and the account of trust developed by Nguyen (2022) is attitudinal and does not need to be directed at people. According to him, trust is about adopting an unquestioning attitude and internalising the reliability of that which is trusted.[21] Such accounts would predict that mathematicians can and indeed do sometimes trust computers.

Setting aside such other accounts of trust, we can ask whether a large-scale collaboration like the Equational Theories Project does in fact embody the (interpersonal) trust-free ideal of mathematics after all. Our reply would be to grant that the role of trust changes, but to insist that it does not disappear. There are still places where mathematicians involved in these projects are relying on others to meet their commitments.

Consider the phase when work gets allocated. Individuals have to register that they will work on particular parts of the proof to avoid duplication. By signing up to some particular set of equational theories, the mathematician makes a commitment to actually work on them, and other collaborators will rely on them to fulfil that commitment to get the overall project finished. More importantly for the integrity of the proof, there is also the use of the Mathlib library, which relies on moderators to check and approve new additions. Implicit in using Mathlib, then, mathematicians must put their trust in the moderators, other contributors to the library, and maintainers of the software.

---

[20] See (Coliva forthcoming) for an explicit discussion of trusting AI-systems.

[21] See (Tacca 2026) for an overview of whether it is possible to trust AI-systems.

## 6. Neural AI: Potentials and Concerns

In sections 4 and 5 we examined the division of epistemic labour among large groups of human mathematicians, facilitated and enabled by technologies such as ITPs. In this section, we consider the philosophical issues that arise when the division of labour in mathematics includes neural AI systems. The basis of these new technologies is not symbolic derivations in a specific formal system, but statistical learning from large datasets, typically via deep neural networks. With the rapid rise of LLMs and deep learning, the ability of these systems to do mathematics is in the spotlight, particularly as one of the major hurdles for the models is to be able to *reason*, and mathematics is often seen as the paradigm test case of abstract reasoning.

Deep learning relevant to mathematics was already being developed to work in more constrained epistemic environments by *DeepMind* in their Go models *AlphaGo* and *AlphaZero*. However, writing good mathematics is a significantly harder challenge due to the large possibility space and the complexity of mathematical language, which typically isn't formal but combines mathematical symbols and notations with logical manipulations and parts of natural language. Modern generative AI systems – large language models (LLMs) and related architectures with billions of parameters – are trained primarily on text and code and are optimized to predict or generate plausible continuations. While originally developed for natural language tasks, they are increasingly being repurposed in mathematics.

LLMs, like ChatGPT, Gemini, and Claude, at the time of writing, seem to be reasonably good at mathematics up to undergraduate level. Numerous articles have attempted to assess the quality of the mathematics outputted by LLMs (e.g. Plevris et al. 2023; Bubeck et al. 2023; Collins et al. 2024). A notable neural AI model is DeepMind's *AlphaProof*, a neural theorem-proving system that integrates large neural networks with reinforcement learning to find proofs. It is an "agent that learns to find formal proofs through RL [reinforcement learning] by training on millions of auto-formalized problems" (Hubert et al. 2025, 1). DeepMind reports (Hubert et al. 2025) that AlphaProof, combined with AlphaGeometry 2 (another of DeepMind's AI agents, focused on solving geometry problems), was able to achieve a result equivalent to a Silver Medal at the 2024 *International Mathematical Olympiad* (IMO). This was a surprising achievement.

And the story does not end there: one year later, at the 2025 IMO, multiple companies claimed that their models generated full answers for five out of six questions, which would be equivalent to a Gold Medal performance.[22] Interestingly, Google *DeepMind* and an experimental model from *OpenAI* both produced natural language solutions, while *Harmonic*'s Aristotle model and *ByteDance*'s Seed Prover model produced formally verifiable *Lean* code solutions. This is certainly an amazing breakthrough, although some mathematicians have been raising concerns about the testing methodology (Riehl 2025).[23] For example, Tao (2025a) explains that, for humans to obtain a Gold Medal, many constraints must be met – the solution must be found within a certain time limit, by a single person, etc. – and it is unclear whether AI models meet what could be considered the equivalents of these constraints.

[22] It should be noted that these results explicitly used a separate model for answering a geometry question.

In more ambitious work, researchers have leveraged neural models to discover new conjectures, patterns, and even proofs, thereby addressing human cognitive limitations and potentially surpassing human mathematicians on specific tasks (Bengio and Malkin 2024, 458). An early result in this direction is the discovery in knot theory of a new knot invariant (Davis 2021).

Romera-Paredes et al. (2024) report using an LLM as part of a system called *FunSearch* to achieve new lower bounds for some versions of the cap set problem. Their idea was to combine the LLM with an evaluator. Then they had the LLM write programs for generating the sets and the evaluator assessing how well the programs do. The best programs were then used as a basis for generating new programs, which were then evaluated again, repeating this process over multiple rounds of an evolutionary algorithm. Crucially, the fact that the LLM is generating programs that are then evaluated, means that there is little danger of AI hallucinations. The role of the LLM is one of generating novel program ideas, as the authors point out: "the LLM should [...] be seen as a source of diverse (syntactically correct) programs with occasionally interesting ideas." (Romera-Paredes et al. 2024, p. 473).

These systems can also support mathematical practice in several ways. First, they can automate routine or tedious tasks: generating LaTeX code, rewriting arguments, and suggesting examples and counterexamples. Second, they can assist with verification and exploration: for instance, by suggesting proof strategies, filling in missing steps in a partially formalized proof, or being used in combination with ITPs and ATPs. Third, they can provide LLM-aided literature search. The mathematical literature is so large that it is unwieldy, and even slight changes in terminology or notation can lead to traditional search engines missing relevant results. Recently, this contributed to the library of Erdös problems, which collects all of the problems attributed to Paul Erdös, along with the status as open or resolved. An

[23] For an excellent discussion of these results, see Kevin Buzzard's blogpost: https://xenaproject.wordpress.com/2025/08/03/ai-at-imo-2025-a-round-up/
Buzzard stresses the fact that the results are unverifiable in the sense that there is no way to guarantee that the companies haven't cheated in some way, presumably to maximise shareholder value. Furthermore, despite press releases suggesting otherwise, there were no particular rules laid down, nor evidence provided of how the results were achieved, and none of the companies have submitted their work to peer review (although the answers were made available). As such, the level of achievement is somewhat opaque, and rigorous citations are woefully lacking.

LLM literature search helped to identify solutions of partial results for multiple problems that had previously been listed as open.[24]

Finally, one envisioned usage for LLMs is autoformalisation, such that informal mathematics can be automatically formalised and then verified by a formal proof checker. For example, Wu et al. (2022) gave an early attempt at autoformalisation of mathematical statements into the formal language of *Isabelle*, and further progress by Patel et al. (2023) and Lu et al. (2024), among others. A major obstacle to general autoformalisation for mathematics is that formal proofs are fragile: one small typo can invalidate a formal proof. A practical solution to this would be to make the formalisation process interactive rather than automated.[25] An example of this in practice is the recent "vibe formalisation" of a proof about Sidon sets by Alexeev and Mixon (2025), meaning they formalised the proof in Lean with the help of an LLM without knowing Lean code themselves. They say:

> It was clear to us that Lean would help us determine the truth, but we didn't know Lean, and it isn't terribly user-friendly. However, ChatGPT can write Lean, so we decided to vibe code the whole proof. It took a long time (about a week), and the process was extremely annoying, but somehow it succeeded. (Alexeev and Mixon 2025, 12)

As a proof of concept, this is an important contribution, as one of the standard reasons against more widespread adoption of formal proofs is the high barrier to entry. If it is possible to accomplish full formalization with *Lean* without knowing the formal language, then making the process of formalizing less annoying makes it a significantly less formidable task.

One set of criticisms of the use of LLMs in mathematics comes from Tanswell and Berg (2026), who survey the philosophical questions that arise from the use of LLMs concerning the status of mathematical knowledge and justification. Their first argument is that LLMs for mathematics have an alignment problem: that when doing mathematics beyond their direct training data, there is no mathematical guarantee that they will use the same conceptual apparatus as human mathematicians. As in Wittgenstein's rule-following problem, the fact that the LLMs are only trained on finitely many instances of seeing a mathematical concept being used means that their statistical predictions of how it will continue to be used may well fail to attach to the obvious way that humans would easily agree on thanks to their shared context and training. Therefore, proofs produced by LLMs cannot be relied upon to be correct beyond their training data.

[24] For a summary, see the Mastodon post by Terence Tao: https://mathstodon.xyz/@tao/115385022005130505

[25] Another sensible approach is to not make autoformalisation a one-shot process, but a modular and iterated one that checks for errors and whether the proof succeeds in the relevant system.

One solution would be to insist on human checking: that once a human has verified the proof it can be relied upon. However, Tanswell and Berg's second argument is that human checking cannot be a systematic defence against the kind of errors that LLMs introduce, for simple psychological reasons. The proofs written by LLMs are produced according to next-token prediction, so sound convincing in a way that can be seductive and authoritative if one is not paying attention to the details. Humans are not good at constant, perfect vigilance, so it is easy for them not to be paying attention. Human mathematicians have learned techniques for reading mathematics to identify errors, but these will be likely ineffective for LLM-produced proofs because the kind of errors that LLMs make are different to the kind of errors that humans make and, crucially, those that mathematicians have been trained to spot effectively. There is also a connection between the two arguments: if the models do seem to be performing reliably beyond their training data, this will lead to human checkers being less thorough and attentive.

Beyond these worries about alignment, there is the more typical worry about hallucinations, or the propensity for LLMs to output authoritative sounding "bullshit" (Hicks, Humphries, and Slater 2024). This is a general problem for LLMs but applies to the mathematics case too. The current models generally continue producing text even if they have made an error or are unable to find a good next step in the proof, meaning that at a practical level they are often unhelpful as proving assistants because they cannot be relied upon to recognise the limits of what they can do, unlike most human collaborators. Of course, this is something that can be addressed with formal verification, as a false proof will not pass it. Again, the most promising prospect at the moment consists of effectively combining symbolic and neural AI systems.

These are worries within mathematics. However, there are also numerous other reasons to be resistant or at least more cautious about the use of AI. First, the environmental damage caused by the massive amounts of power, water and space needed for the training and use of the models is morally unacceptable when the world is facing climate catastrophe. Second, there is a danger of mathematicians becoming over-reliant on AI systems, and thereby collectively deskilling mathematics. The productive struggle of being stuck in mathematics may not be a hindrance to mathematical thinking, but a crucial part of it that, when removed, leads to a degeneration of mathematical abilities. Third, many are worried about the training data used to develop the large commercial LLMs and its provenance, such as the use of pirated copyright materials without acknowledgement or remuneration. Likewise, the mathematical training likely involved available online sources like the ArXiv, MathOverflow, and the various formal proof libraries. These were produced by mathematicians as part of the collective endeavour of mathematics, but are used to make a profit by big technology companies. Fourth, and relatedly, the values of the tech industry are not aligned with those of professional mathematicians, so the ceding of control over the future of mathematics endangers the values that mathematicians hold (as is argued by Harris (2024)). Fifth, another worry is about the cybersecurity of LLMs, which are widely known to be vulnerable to adversarial attacks. Formal mathematics could be intentionally "hacked." If and when these tools are integrated more deeply into mathematics, this could create unanticipated cybersecurity problems, such as for the integrity of the mathematical literature or the Lean Mathlib. Sixth, this also extends to the problem of it becoming easier to generate plausible-seeming mathematical spam that could clog up existing systems of mathematical peer-review and dissemination, such as the ArXiv.

To summarise, the use of AI in mathematics is developing quickly in a variety of ways: to automatically generate mathematics, to formalise and verify proofs, and as an interactive assistant. Nonetheless, there are also serious concerns that suggest that caution is necessary in dealing with these new tools.

## 7. Conclusion

The technological turn in mathematics involves a dual transformation, affecting both the epistemology and social organization of mathematics. Interactive Theorem Provers, symbolic automated reasoning systems, and neural AI tools do not merely accelerate existing practices; they reshape how proofs are produced, checked, and distributed. The emerging interplay between these different technologies, and in particular between symbolic and neural approaches to AI, suggests that no single technological paradigm will dominate mathematical practice; rather, future mathematics is likely to be shaped by complex combinations of formal verification, automation, and statistical learning.

Many of the concerns these technologies raise are continuous with long-standing issues in the philosophy of mathematics. For instance, they pose questions about the reliability of software, the fallibility of mathematical reasoning, and the nature of collaborations. Still, the new ways in which ITPs are used and the opacity of contemporary neural AI systems introduce genuinely novel challenges. Moreover, these technologies increasingly function not merely as passive instruments but as active participants in mathematical practice, redistributing epistemic labour between humans and machines. This raises pressing questions about agency, trust, and responsibility, and about authorship.[26] Can mathematical knowledge be appropriately ascribed to individuals, collectives, and hybrid human–machine systems?

These developments have implications not only for the philosophy of mathematics, but also for social epistemology and the philosophy of science more broadly. Understanding this evolving landscape is essential if we are to assess both the promises and the risks of AI-driven mathematics, and to shape its future in ways that remain aligned with the values of mathematical inquiry.

[26] See, for example, (Habgood-Coote et al. 2024).

## Acknowledgements

We wish to thank Ásgeir Berg, Jacopo Emmenegger, Axel Gelfert, Ursula Martin, Alexander Paseau, and Tito Sacchi for their helpful feedback. Part of this research was funded by the Italian Grant FIS2, *Humanizing Mathematical Knowledge: Fallibility, Technology, Know How* (FIS-04053, PI: De Toffoli) - CUP I53C24003170001.